\documentclass[a4paper]{article}

\usepackage[english]{babel}
\usepackage[utf8]{inputenc}

\usepackage{tikz-cd}
\usepackage{tikz}
\usetikzlibrary{matrix,arrows}
\usepackage[colorinlistoftodos]{todonotes}

\usepackage{amscd, amsmath, amssymb, amsthm, enumerate, mathrsfs, graphicx, stackengine, amsfonts}

\usepackage{relsize}
\usepackage[all,cmtip]{xy}
\usepackage[margin=1.5in]{geometry}
\usepackage[normalem]{ulem}
\usepackage[utf8]{inputenc}
\usepackage[T1]{fontenc}
\usepackage{lmodern}
\usepackage{lmodern}
\usepackage{lipsum}

\let\oldabstract\abstract
\let\oldendabstract\endabstract
\makeatletter
\renewenvironment{abstract}
{%
               {\list{}{\addtolength{\leftmargin}{1em} 
                        \listparindent 1.5em%
                        \itemindent    \listparindent%
                        \rightmargin   \leftmargin%
                        \parsep        \z@ \@plus\p@}%
                \item\relax}%
               {\endlist}%
\oldabstract}
{\oldendabstract}
\makeatother

\newcommand{\RN}[1]{%
  \textup{\uppercase\expandafter{\romannumeral#1}}%
}

\title{Weak Bounded Negativity Conjecture}

\author{Feng Hao}
\date{\vspace{-5ex}}
\begin{document}
\maketitle
\begin{abstract} 

\noindent In this paper, we prove the following ``Weak Bounded Negativity Conjecture'', which says that given a complex smooth projective  surface $X$, for any reduced curve $C$ in $X$ and integer $g$, assume that the geometric genus of each component of $C$ is bounded from above by $g$, then the self-intersection number $C^2$ is bounded from below.
\end{abstract}

\section{Introduction}

The so called Weak Bounded Negativity Conjecture (Conjecture 1.2) is  motivated by the study of the old folklore conjecture ``Bounded Negativity Conjecture'', which is stated as follows. 

\textbf{Conjecture 1.1 (Bounded Negativity Conjecture):} \textit{For any smooth complex projective surface $X$, there exists a constant $b(X)$ only depending on $X$ itself, such that $C^2\geq b(X)$ for any reduced curve $C$ in $X$.}\\

In this paper, we consider the following Weak Bounded Negativity Conjecture.

\textbf{Conjecture 1.2 (Weak Bounded Negativity Conjecture):}  \textit{For any smooth complex projective surface $X$ and any integer $g$, there is a constant $b(X,g)$ only depending on $X$ and integer $g$, such that $C^2\geq b(X, g)$ for any reduced curve $C=\Sigma C_i$ in $X$ with the geometric genus $g(C_i)\leq g$, for all $i$.}\\

For the Weak Bounded Negativity Conjecture, there are several partial results as stated in Theorem 1.3 and Theorem 1.4.

\textbf{Theorem 1.3 (Bogomolov):}  \textit{Let $X$ be a smooth projective surface with Kodaira dimension $\kappa(X)\geq 0$. Then for any smooth irreducible curve $C\subset X$ of geometric genus $g(C)$,  we have\[C^2\geq K_X^2-4c_2(X)-4g(C)+4,\]
where $c_1$ and $c_2$ are the first and second Chern numbers of the surface $X$, respectively.}

The proof of Theorem 1.3 involves Bogomolov's criterion for the unstable bundles on surfaces and Bogomolov-Sommese vanishing Theorem. Refer to \textbf{[Bau2, Theorem 3.4.4]} and Bogomolov \textbf{[Bogo, section 5]} for details.\\

Th. Bauer, B. Harbourne,  T. Szemberg, and other authors used the logarithmic Miyaoka-Yau inequality to prove the following theorem and gave a better bound.

\textbf{Theorem 1.4 (\textbf{[Bau2, Theorem 2.6]}):}  \textit{Let $X$ be a smooth projective surface with Kodaira dimension $\kappa(X)\geq 0$. Then for any integral curve $C\subset X$,  we have\[C^2\geq K_X^2-3c_2(X)+2-2g(C)\]
where $c_1$ and $c_2$ are the first and second Chern numbers of surface $X$, respectively.}\\

\textbf{Remark 1.5:} As T. Szemberg mentioned to me, Theorem 1.4 is actually a corollary of the generalized Logarithmic Miyaoka-Yau Inequality: Miyaoka \textbf{[Miy, Theorem 1.1]}. We will use \textbf{[Miy, Theorem 1.1]} in section 2.\\

In this paper, we use the elementary intersection theory, the generalized Logarithmic Miyaoka-Yau Inequality (\textbf{[Miy, Theorem 1.1]}), and some techniques in the proof of Theorem 1.4  to give a full proof of the Weak Bounded Negativity Conjecture.

\section{Integral curves in any smooth complex projective surface $X$}

In this section, we prove the Weak Bounded Negativity Conjecture for integral curves in a surface $X$ through case by case analysis, looking at  $H^0(X, \mathcal{O}_X(-K_X))$, where $\mathcal{O}_X(-K_X)$ is the anti-canonical line bundle of $X$. We denote $dim H^0(X, \mathcal{O}_X(D))$ by $h^0(D)$, for any divisor $D$ on $X$.

\subsection{ Surface $X$ with $h^0(-K_X)>0$}

For this case, we have the following simple observation.\\

\textbf{Lemma 2.1.1:}\ \ \textit{Given a smooth projective surface $X$ over $\mathbb{C}$ with $h^0(-K_X)>0$, then the Weak Bounded Negativity Conjecture holds for integral curves in $X$.}

\textit{Proof.} Since $h^0(-K_X)>0$, we can choose an effective divisor in the linear system $|-K_X|$ and still call it $-K_X$. Note that it contains only finitely many integral negative curves with negative self-intersection. Then for any other integral curve $C$ which are not in the components of $-K_X$, we have by genus formula$$g_a(C)=1+1/2(C^2+C\cdot K_X),$$ where $g_a(C)$ is the arithmetic genus of C. 
Therefore we have $C^2=2g_a(C)-2-C\cdot K_X$. Since $-K_X$ is effective and $C$ is not in the components of $-K_X$, we have $C^2\geq -2$. Note that the bound in this case does not depend on the geometric genus of $C$.

$\blacksquare$ \\

\textbf{Example 2.1.2:}\ Consider the minimal rational surfaces: Hirzebruch surfaces $\Sigma_n$.

Note first $K_{\Sigma_n}^2=8$. By Riemann-Roch formula, we have $$h^0(-K_{\Sigma_n})+h^0(2K_{\Sigma_n})-h^1(-K_{\Sigma_n})=1+K^2_{\Sigma_n}.$$ Since $\Sigma_n$ is a rational surface, we have $h^0(2K_{\Sigma_n})=0$. Thus we get $h^0(-K_{\Sigma_n})\geq 9$. On the other hand, we know that there is only one negative curve on $\Sigma_n$, with self-intersection $-n$. Then by the above lemma, we know that if $n>2$, the negative curve is contained in an effective representative of $-K_{\Sigma_n}$.\\

\subsection{Surface $X$ with $h^0(-K_X)=0$}

In this subsection, we will use the invariant $h^0(m(K_X+C))$ of curves on a surface to divide the problem into two cases.\\

\textbf{Case \RN{1}: $C$ is an integral curve on $X$, such that $h^0(m(K_X+C))=0$ for all $m$.}

\textbf{Lemma 2.2.1:} \textit{ Given a smooth projective surface $X$ with $h^0(-K_X)=0$, and an integral curve $C\subset X$ of arithmetic genus $g_a(C)$ with $h^0(m(K_X+C))=0$ for all $m$, we have $C^2\geq K_X^2+\chi(\mathcal{O}_X)-3$.}

\textit{Proof.} In this case,  $h^0(2(K_X+C))=0$. Hence we have $h^0(2K_X+C)=0$. By Riemann-Roch formula and genus formula for curves on surfaces, we have $$h^0(2K_X+C)+h^0(-K_X-C)-h^1(2K_X+C)=K_X^2+3g_a(C)+\chi(\mathcal{O}_X)-3-C^2.$$

Since $h^0(-K_X)=0$ and $C$ is effective, $h^0(-K_X-C)=0$. Thus we get  $$C^2\geq K_X^2+3g_a(C)+\chi(\mathcal{O}_X)-3\geq K_X^2+\chi(\mathcal{O}_X)-3.$$

$\blacksquare$ \\

\textbf{Remark 2.2.2:} The lower bound in the above case does not involve the geometric genus of the integral curves.\\

\textbf{Case \RN{2}: $C$ is a smooth irreducible curve on $X$, such that $h^0(m(K_X+C))\neq 0$ for some $m>0$.}

For this case, we first introduce the following two theorems:

\textbf{Theorem 2.2.3(Zariski Decomposition Theorem):} \textit{Let $X$ be a smooth projective surface, and let $D$ be a pseudo-effective integral divisor on $X$. Then $D$ can be written uniquely as a sum
$D=P+N$ of $\mathbb{Q}$-divisors with the following properties:}

\textit{(1) $P$ is nef;}

\textit{(2) $N=\Sigma^r_{i=1}a_iE_i$ is effective, and if $N\neq 0$ then the intersection matrix
$$\|E_i\cdot E_j\|$$
determined by the components of N is negative definite;}

\textit{(3) $P$ is orthogonal to each of the components of $N$ , i.e.  $P\cdot E_i=0$.}\\

Refer to Fujita \textbf{[Fuj, Theorem 1.12]} for the proof of Theorem 2.2.3.\\

\textbf{Remark 2.2.4:} The above version of the Zariski Decomposition Theorem is due to Fujita. The original version of the Zariski Decomposition Theorem says that for an effective divisor $D$,  we have a unique decomposition $D=P+N$ satisfying the above three properties, where $P$ is necessarily effective. Refer to Zariski \textbf{[Zar, Theorem 7.7]} for the original version. \\

\textbf{Theorem 2.2.5:}  \textit{Let $X$ be a smooth projective surface and $C$ be a smooth curve in $X$. Assume that $K_X+C$ is pseudo-effective. According to Theorem 2.2.3, $K_X+C$ admits a Zariski decomposition. Then the following inequality  holds
$$c_2(X)-e(C)-\frac{1}{3}(K_X+C)^2+ \frac{1}{12}N^2\geq 0,$$
where $e(C)$ is the topological Euler characteristic class of $C$, and $N$ is the negative part (non-nef part) of the Zariski decomposition of $K_X+C$.}\\

Theorem 2.2.5 is a special case of Miyaoka \textbf{[Miy, Theorem 1.1]}. \\


\textbf{Corollary 2.2.6:} \textit{Let $X$ be a smooth projective surface with $H^0(X, -K_X)=0$, and $C\subset X$ be a smooth irreducible curve of genus $g(C)$ with $H^0(X, m(K_X+C))\neq 0$ for some $m$. Then $C^2\geq K_X^2-3c_2(X)+2-2g(C)$.}

\textit{Proof.} Since there exists $m>0$ such that $h^0(m(K_X+C))>0$, $K_X+C$ is a pseudo-effective divisor. By Theorem 2.2.3, $K_X+C$ admits a Zariski decomposition $K_X+C=P+N$, with $P$ the nef part. Then by Theorem 2.2.5, we get the following inequality
$$c_2(X)-e(C)-\frac{1}{3}(K_X+C)^2+ \frac{1}{12}N^2\geq 0.$$ Note that $N^2\leq 0$ by property (2) of the Zariski Decomposition Theorem. Thus we have $(K_X+C)^2\leq 3(c_2(X)-2+2g(C))$. Note also that by genus formula, we have $$(K_X+C)^2=K_X^2+4(g(C)-1)-C^2,$$

Hence $$C^2\geq K_X^2-3c_2(X)+2-2g(C).$$ \ $\blacksquare$ \\

Next we will modify the strategy in the proof of Theorem 1.4 to prove the Weak Bounded Negativity Conjecture for integral curves. Considering Lemma 2.2.1 and Corollary 2.2.6, to prove the Weak Bounded Negativity Conjecture for integral curves, it suffices to prove the following theorem.

\textbf{Theorem 2.2.7:} \textit{Let $X$ be a smooth projective surface with $H^0(X, -K_X)=0$, and $C\subset X$ be an integral curve with $H^0(X, m(K_X+C))\neq 0$ for some $m$. Then $C^2\geq \mu(X, g(C))$, where $\mu(X, g(C))$ is defined to be $min\{K_X^2+\chi(\mathcal{O}_X)-3,\   K_X^2-3c_2(X)+2-2g(C)\}$.}\\

First we have the following simple observation.

\textbf{Lemma 2.2.8:} \textit{Let $X$ be a smooth projective surface with $H^0(X, -K_X)=0$, and $p$ be a point in $X$. Let $\pi: \tilde{X}=Bl_p(X)\rightarrow X$ be the blow up along $p$. Then $H^0(X, -K_{\tilde{X}})=0$.}

\textit{Proof.} Since $h^0(-K_X)=0$, $h^0(-\pi^*(K_X))= 0$. Note that $-K_{\tilde{X}}=-\pi^*K_X-E$, where $E$ is the exceptional divisor of the blow up.  Thus $h^0(-K_{\tilde{X}})=0$.\ $\blacksquare$ \\

By Lemma 2.2.8, given a surface $X$  with $h^0(-K_X)=0$, a blow up of $X$ will satisfy the same property. Thus to prove theorem 2.2.7, it suffices to prove the following lemma.

\textbf{Lemma 2.2.9:} \textit{Let $X$ be a smooth projective surface with $h^0(-K_X)=0$, $C$ be an integral curve of geometric genus $g(C)$, and $p\in C$ be a point with $m:=mult_pC\geq 2$. Let $\pi: \tilde{X}\rightarrow X$ be the blow up of X at $p$ with the exceptional divisor $E$. Let $\tilde{C}=\pi^*(C)-mE$
be the strict transform of $C$. Then the inequality $$\tilde{C}^2\geq \mu(\tilde{X}, g(\tilde{C}))$$
implies $$C^2\geq \mu(X, g(C))$$}

\textit{Proof.} Note that we have $C^2=\tilde{C}^2+m^2,\  K^2_X=K^2_{\tilde{X}}+1, \ c_2(X)=c_2(\tilde{X})-1, \ \chi(\mathcal{O}_X)=\chi(\mathcal{O}_{\tilde{X}})$, and\  $g(C)=g(\tilde{C}).$

Note that $K_X^2+\chi(\mathcal{O}_X)-3$ and  $K_X^2-3c_2(X)+2-2g(C)$ only depend on $X$ in the blow-up procedure. Thus we may denote $$M(X)=K_X^2+\chi(\mathcal{O}_X)-3$$ and $$N(X)=K_X^2-3c_2(X)+2-2g(C).$$ Then we have $$M(\tilde{X})+1=M(X)$$ and $$N(\tilde{X})+4=N(X).$$

There are three cases

(1) If $\mu(X, g(C))=M(X)$ and $\mu(\tilde{X}, g(\tilde{C}))=M(\tilde{X})$, \ \ $C^2-m^2=\tilde{C}^2\geq M(\tilde{X})=M(X)-1$ implies $C^2\geq M(X)$.

(2) If $\mu(X, g(C))=M(X)$ and $\mu(\tilde{X}, g(\tilde{C}))=N(\tilde{X})$, \ \ $C^2-m^2=\tilde{C}^2\geq N(\tilde{X})=N(X)-4$ implies $C^2\geq N(X)\geq M(X)$.

(3) If $\mu(X, g(C))=N(X)$ and $\mu(\tilde{X}, g(\tilde{C}))=N(\tilde{X})$, \ \ $C^2-m^2=\tilde{C}^2\geq N(\tilde{X})=N(X)-4$ implies $C^2\geq N(X)$. \ $\blacksquare$ \\

\textbf{Proof of Theorem 2.2.7:} It follows immediately from Lemma 2.2.9. \ $\blacksquare$ \\

\textbf{Corollary 2.2.10:} \textit{The Weak Bounded Negativity Conjecture holds for integral curves.}

\textit{Proof.} Collect the results: Lemma 2.1.1, Lemma 2.2.1, and Theorem 2.2.7, we get the above corollary. \ $\blacksquare$ \\

\section{Reduced curves with arbitrary singularity}

In this section, we will prove the Bounded Negativity Conjecture for the general reduced curves case, based on the results we get in section 2. However, the idea of the proof of the following theorem comes from the proof of \textbf{[Bau1, Theorem 5.1]} with  a different situation. 

\textbf{Theorem 3.1.1:} \textit{Let $C=\Sigma_i C_i$ be a reduced curve on any smooth complex projective surface $X$. Suppose that there exists an integer $g$ such that the geometric genus $g(C_i)\leq g$ for all $i$. Then there exists a constant $B(X, g)$ only depending on $X$ and $g$, such that $$C^2\geq B(X,g).$$}

\textit{Proof.} In Remark 2.2.4, we have the Zariski Decomposition Theorem (\textbf{[Zar, Theorem 7.7]}). Then $C=P+N$, where $P$ is nef and effective and $N=\Sigma^r_{i=1}a_iE_i$ is effective. \ Then $$C^2=(P+\mathlarger{\mathlarger{\sum}}_{i=1}^{r}a_iE_i)^2=P^2+(\mathlarger{\mathlarger{\sum}}_{i=1}^{r}a_iE_i)^2\geq (\mathlarger{\mathlarger{\sum}}_{i=1}^{r}a_iE_i)^2.$$

Since $C$ is reduced and $P$, $N$ are effective, we have $a_i\leq 1$. By Hodge Index Theorem, and matrix $[E_i\cdot E_j]$ is negative definite, we have $r\leq h^{1,1}-1$. Also, by Corollary 3.0.4, $E_i^2\geq b(X, g)$ for some constant $b(X,g)$ depending on $X$ and $g$. Also, we can always assume $b(X,g)\leq0$.

Thus we get $C^2\geq a_1^2E_1^2+...+a_r^2E_r^2\geq (h^{1,1}(X)-1)b(X,g).$ Then just let $B(X, g)=(h^{1,1}(X)-1)b(X,g).$

\ $\blacksquare$ \\








\newpage

\textbf{References}

\textbf{[Bau1]} Bauer, Th., et al.  Negative curves on algebraic surfaces. arXiv.org/abs/1109.1881.

\textbf{[Bau2]} Bauer, Th., et al. Recent developments and open problems in linear series. To appear in “Contributions to Algebraic Geometry”, Impanga Lecture Notes Series. arXiv:1101.4363.


\textbf{[Bogo]} F. A. Bogomolov. Stable vector bundles on projective surface. Russian Academy of Sciences. Sbornik Mathematics.  Volume 81, Number 2.

\textbf{[Fuj]} T. Fujita. On Zariski Problem. Proc. Japan Acad., 5, Ser. A (1979) Vol. 55(A).

\textbf{[Miy]} Y. Miyaoka. The Maximal Number of Quotient Singularities on Surfaces with Given Numerical Invariants. Mathematische Annalen (1984)
Volume: 268, page 159-172.






\textbf{[Zar]} O. Zariski, The Theorem of Riemann-Roch for High Multiples of an Effective Divisor on an Algebraic Surface. Ann. of Math, 1962




\newpage

Feng Hao

Department of Mathematics, Purdue University

150 N. University Street, West Lafayette, IN 47907-2067

E-mail address: \textit{fhao@purdue.edu}

\end{document}